\documentclass[submission]{FPSAC2026}

\usepackage{graphics}
\usepackage{mathrsfs}
\usepackage{tikz-cd}
\usepackage{mathabx}
\usepackage{comment}

\usepackage{amsthm}
\usepackage{aliascnt}
\usepackage{cleveref}

\newtheorem{theorem}{Theorem}[section]

\newcommand{\newtheoremwithalias}[2]{%
  \newaliascnt{#1}{theorem}%
  \newtheorem{#1}[#1]{#2}%
  \aliascntresetthe{#1}%
  \crefname{#1}{#2}{#2s}%
  \Crefname{#1}{#2}{#2s}%
}

\newtheoremwithalias{lemma}{Lemma}
\newtheoremwithalias{corollary}{Corollary}
\newtheoremwithalias{problem}{Problem}
\newtheoremwithalias{conjecture}{Conjecture}

\theoremstyle{definition}
\newtheoremwithalias{definition}{Definition}

\newtheoremwithalias{example}{Example}

\newcommand{\wt}[1]{\widetilde{#1}}

\newcommand{\Z}{\mathbb{Z}}

\newcommand{\A}{\mathcal{A}}
\newcommand{\B}{\mathcal{B}}
\newcommand{\C}{\mathcal{C}}
\renewcommand{\P}{\mathcal{P}}

\newcommand{\cR}{\mathcal{R}}
\newcommand{\id}{\mathrm{id}}

\newcommand{\cG}{\mathcal{G}}
\newcommand{\CR}{\mathrm{Cr}}
\renewcommand{\S}{\mathfrak{S}}
\newcommand{\affineS}{\wt{\S}}

\newcommand{\defn}[1]{{\emph{#1}}}

\DeclareMathOperator{\Inv}{Inv}

\DeclareMathOperator{\Rev}{Rev}

\usepackage{lipsum}

\title[Commutation classes of reduced words and higher Bruhat orders for affine permutations]{Commutation classes of reduced words and \\ higher Bruhat orders for affine permutations}

\author{Sara Billey\thanks{\href{mailto:billey@uw.edu}{billey@uw.edu}}\addressmark{1}, Herman Chau\thanks{\href{mailto:herman.h.chau@gmail.com}{herman.h.chau@gmail.com}}\addressmark{1}, \and Kevin Liu\thanks{\href{mailto:keliu@sewanee.edu}{keliu@sewanee.edu}}\addressmark{2}}

\address{\addressmark{1}Department of Mathematics, University of Washington, Seattle, WA, USA\\ \addressmark{2}Department of Mathematics and CS, The University of the South, Sewanee, TN, USA}

\received{November 10, 2025}

\revised{March 9, 2026}

\abstract{The higher Bruhat orders are partial orders that generalize the weak order on the symmetric group $\mathfrak{S}_n$, and the second higher Bruhat order is a poset on commutation classes of reduced words for the longest element in ${\mathfrak{S}}_n$, where covering relations correspond to braid relations. Constructing analogs in other settings is an area of recent interest, and we present an analog that generalizes any interval $[\id,w]$ in the weak order of the affine symmetric group $\widetilde{\mathfrak{S}}_n$. Paralleling the classical case, we show the second higher Bruhat order is a poset on commutation classes of reduced words for $w\in \widetilde{\mathfrak{S}}_n$. When $w\in \mathfrak{S}_n$, we also establish results for all higher Bruhat orders that are direct analogs of ones in the classical case. }


\usepackage[backend=bibtex]{biblatex}
\begin{document}

\maketitle

\section{Introduction}\label{sec:intro}

The classical {higher Bruhat orders} $\B(n,k)$ are partial orders introduced by Manin
and Schechtman \cite{manin_schechtman_1989} in their study of fundamental groups of certain hyperplane arrangements. They were further characterized by Ziegler \cite{ziegler_1993} {using single step inclusion of consistent sets}.  The first higher Bruhat order $\B(n,1)$ is isomorphic to the weak (Bruhat) order on the
symmetric group $\S_n$, and the second higher Bruhat order $\B(n,2)$ is
a partial order on commutation classes of reduced words for the
longest permutation.  In the general case, $\B(n,k)$ can be described using pseudo-hyperplane arrangements
\cite{ziegler_1993}, oriented matroids \cite{sturmfels-oriented}, and
zonotopal tilings \cite{felsner-zonotopes}. The higher Bruhat orders have applications to Bott-Samelson varieties
\cite{Elias_2016} and Steenrod algebras \cite{steenrod}.  

Generalizing
the higher Bruhat orders to other settings has been an area of recent
interest \cite{danilov2022BC,higherbruhatB}. In
2022, Ben Elias generalized the second higher Bruhat orders to include
nonreduced words to prove a
generalized Bergman diamond lemma for Hecke-type algebras
\cite{Elias.2022}.  Daniel Hothem generalized the
higher Bruhat orders beginning from intervals
 $[\mathrm{id},w]$ in the weak order on $\S_n$
\cite{hothem2021extendinghigherbruhatorders}.  Elias and Hothem posed
the following question.

\begin{problem}\label{prob:EliasHothem}(Elias-Hothem, 2020)
    Generalize the higher Bruhat orders to all intervals $[\mathrm{id},w]$ in the affine symmetric groups $\affineS_n$.
\end{problem}

For each affine permutation $w\in \affineS_n$, Elias proposed orienting edges between commutation classes of reduced words for $w$ according to the directed braid relations $s_is_{i+1}s_i\to s_{i+1}s_is_{i+1}$, where each index is modulo $n$. This results in a directed graph, denoted $\cG(w)$. 

\begin{conjecture}(Elias, 2021)\label{conj:Elias.n=2}
    For any $w\in \affineS_n$, the directed graph $\cG(w)$ is acyclic with a unique source vertex and a unique sink vertex.
\end{conjecture}

In this extended abstract, we prove \cref{conj:Elias.n=2} and propose a construction for \cref{prob:EliasHothem}. We begin with the set $\Inv_k(w)$ of $k$-inversions for $w\in \affineS_n$. These are the length $k$ decreasing subsequences in the one-line notation of $w$, up to a congruence shift. Paralleling the classical constructions in \cite{manin_schechtman_1989}, we define the set of admissible orders $\A_w(n,k)$ to be certain total orders of $\Inv_k(w)$, introduce a notion of commutation equivalence $\sim_w$, and construct the higher Bruhat orders $\B_w(n,k)$ as a poset on $\A_w(n,k)/\!\sim_w$. We associate to each admissible
order a {reversal set} that is constant up to commutation equivalence. These satisfy conditions that motivate a notion of consistent sets that extends the one due to Ziegler \cite{ziegler_1993}, and we define a poset $\C_w(n,k)$ on the consistent subsets of $\Inv_{k}(w)$ ordered by single step inclusion. From our construction, $\B_w(n,1)\cong \C_w(n,2)$ is isomorphic to the interval $[\mathrm{id},w]$ in the weak order on $\affineS_n$, and we establish the following theorem when $k=2$. \cref{conj:Elias.n=2} follows immediately from this result.

\begin{theorem}\label{thm:affine}
    For any positive integer $n$ and $w\in \affineS_n$, the following hold.
    \begin{enumerate}[label=(\alph*)]
    \item There are natural bijections between maximal chains of $\C_w(n,2)$, reduced words for $w$, reflection orders for $w$, and $\A_w(n,2)$.
    \item $\B_w(n,2)\cong \C_w(n,3)$ is a ranked poset with  unique minimal and maximal elements corresponding with reversal sets $\emptyset$ and $\Inv_3(w)$. 
    \item The Hasse diagram of $\B_w(n,2)\cong \C_w(n,3)$ is isomorphic to $\cG(w)$ as a directed graph. Furthermore, the diameter of $\cG(w)$ as an undirected graph is $|\Inv_3(w)|$.
    \end{enumerate} 
\end{theorem}

\cref{thm:affine} has connections to other work on graphs of reduced words. When $w\in \S_n$, a comparison of our constructions shows $\B_w(n,2)$ is the quotient of a ranked poset on all reduced words for $w$ defined by Assaf \cite{Assaf.2019}. The diameter statement in \cref{thm:affine}(c) also extends a result of Gutierres, Mamede, and Santos \cite{GUTIERRES} from $\S_n$ to $\affineS_n$. 

When $w\in \S_n$, we also extend portions of \cref{thm:affine} to arbitrary $k$, resulting in the following theorem that recovers numerous properties of the classical higher Bruhat orders in \cite{manin_schechtman_1989} and \cite{ziegler_1993}. We conjecture that this theorem holds for arbitrary $w \in \affineS_n$. For $1\leq k\leq n\leq 8$, we have computationally verified this for all $w\in \affineS_n$ of length up to and including $15$. For $n=6$, we have further verified through length $21$, and for $n\leq 5$, we have further verified through length $27$.

 \begin{theorem}\label{thm:nonlongest}
    Let $k,n$ be positive integers and $w\in \S_n$. The following
    hold.
    \begin{enumerate}[label=(\alph*)]
    \item $\B_w(n,k)$ is isomorphic as a poset to $\C_w(n,k+1)$, and the isomorphism sends an equivalence class of admissible orders to the reversal set of the class.
    \item $\B_w(n,k) \cong \C_w(n,k+1)$ is a ranked poset with 
    unique minimal and maximal elements corresponding with reversal
    sets  $\emptyset$ and  $\Inv_{k+1}(w)$.  
    \item There is a natural bijection between maximal chains of $\C_w(n,k+1)$ and $\A_w(n,k+1)$.
    \end{enumerate} 
\end{theorem}

 We begin in \cref{sec:prelims} by giving
preliminary information on the affine symmetric group. As our constructions in $\S_n$ are simpler to understand, we will start in \cref{sec:non-longest} by constructing our higher
Bruhat orders for $w\in \S_n$ and outlining our approach for proving
\cref{thm:nonlongest}.  We then move to the affine higher
Bruhat orders for $w\in \affineS_n$ in
\Cref{sec:affine} and outline our approach for \cref{thm:affine}. We also state \cref{conj:affineG_R}, which is the missing link in a full generalization of \cref{thm:nonlongest} to all affine permutations.

\section{Preliminaries}\label{sec:prelims}

In this section, we give preliminaries on the affine symmetric group. See \cite[Section 8.3]{BB} for additional background. Throughout, $n$ is a positive integer.

A bijection $w: \Z \to \Z$ is called \defn{$n$-periodic} if for all $x
\in \Z$, we have $w(x+n) = w(x)+n$. Such a 
bijection is uniquely determined by its values on $[n]=\{1,2,\ldots,n\}$.  The
\defn{affine symmetric group} $\affineS_n$ is the group of all
$n$-periodic bijections $w: \Z \to \Z$ that satisfy $\sum_{i \in [n]}w(i) = \binom{n+1}{2}$, where the group
operation is function composition. An element in $\affineS_n$ is called an \defn{affine permutation}, and in our work, we will express it in \defn{window notation} $(w(1),w(2),\ldots,w(n))$. Note that $\affineS_n$ contains the usual symmetric group $\S_n$ as the subgroup of permutations such that $w([n])=[n]$.

When $n\geq 2$, the group $\affineS_n$ is generated by the \defn{simple transpositions} $s_0,s_1, \ldots, s_{n-1}$ where $s_i$ maps $s_i(i) = i+1$, $s_i(i+1) = i$, and fixes all $j\in \Z$ such that $j \not\equiv i, i+1 \pmod{n}$. We consider the indices of the simple transpositions modulo $n$, so $s_0$ and $s_n$ denote the same simple transposition. 
The simple transpositions in $\affineS_n$ satisfy the following set of minimal relations: (1) $s_i^2 = 1$ for all $i \in \Z$, (2) $s_is_{i+1}s_i = s_{i+1}s_is_{i+1}$ for all $i \in \Z$, and (3) $s_is_j = s_js_i$ for all $i,j \in \Z$ with $j \not\equiv i,i+1, i-1 \pmod{n}$.
Relations of the second kind are termed \defn{braid relations}, and relations of the third kind are termed \defn{commutation relations}.
The symmetric group $\S_n$ is generated by $s_1, \ldots, s_{n-1}$, omitting $s_0$.

A \defn{reduced expression} for an element $w \in \affineS_n$ is a
product $s_{i_1} \cdots s_{i_\ell}$ of minimal length such that
$s_{i_1} \cdots s_{i_\ell} = w$. For brevity, we will consider
\defn{reduced words}, which are the indices $i_1 \cdots i_\ell \in (\Z/n\Z)$ of a reduced expression $s_{i_1} \cdots s_{i_\ell}$. All reduced words for an element $w \in \affineS_n$ have the same \defn{length}, which is denoted $\ell(w)$. 
In general, an element $w\in \affineS_n$ can have many reduced words, and given a reduced word, one can generate additional ones using relations among the generators of $\affineS_n$. 
Two reduced words are \defn{commutation equivalent} if they differ by a sequence of commutation relations. If a reduced word for $w$ has $i(i+1)i$ (resp. $(i+1)i(i+1)$) in consecutive entries, then a \defn{braid} is the operation that replaces these three elements with $(i+1)i(i+1)$ (resp. $i(i+1)i$).

\begin{figure}
    \centering
    \scalebox{0.75}{
    \begin{tikzpicture}
        \node (1) at (0,-3) {$4121432$};
        \node (2) at (2,-3) {$4212432$};
        \node (3) at (4,-3) {$4214232$};
        \node (4) at (6,-3) {$4214323$};
        \node (5) at (2,-1.5) {$2412432$};
        \node (6) at (4,-1.5) {$2414232$};
        \node (7) at (6,-1.5) {$2414323$};
        \node (8) at (4,0) {$2141232$};
        \node (9) at (6,0) {$2141323$};
        \node (10) at (8,0) {$2143123$};
        \path (1) edge node {} (2);
        \path (2) edge node {} (3);
        \path (3) edge node {} (4);
        \path (2) edge node {} (5);
        \path (3) edge node {} (6);
        \path (4) edge node {} (7);
        \path (5) edge node {} (6);
        \path (6) edge node {} (7);
        \path (6) edge node {} (8);
        \path (7) edge node {} (9);
        \path (8) edge node {} (9);
        \path (9) edge node {} (10);
    \end{tikzpicture} \quad 
    \begin{tikzpicture}[every text node part/.style={align=center}]
        \node (1) at (0,-2) {$4121432$};
        \node (2) at (3,-2) {$2414232$ \\ $2412432$ \\ $4212432$ \\ $4214232$  };
        \node (3) at (6,-2) {$2414323$ \\ $4214323$  };
        \node (4) at (3,0) {$2141232$};
        \node (5) at (6,0) {$2141323$ \\ $2143123$};
        \draw[->] (1) edge node {} (2);
        \draw[->] (2) edge node {} (3);
        \path[->] (2) edge node {} (4);
        \path[->] (3) edge node {} (5);
        \path[->] (4) edge node {} (5);
    \end{tikzpicture}
    }
    \caption{$\cR(w)$ (left) and $\cG(w)$ (right) for the affine permutation $w= (1,7,2,0) \in \affineS_4$.}
    \label{fig:reducedwordgraph}
\end{figure}
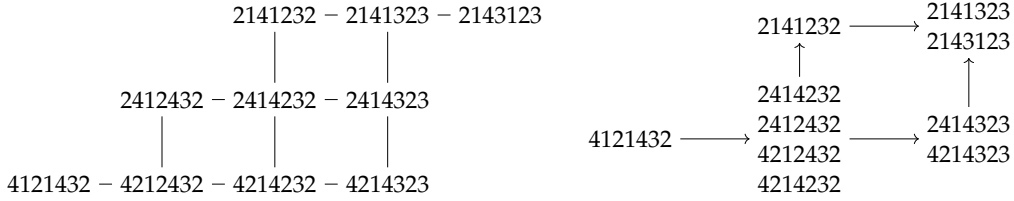

The \defn{graph on reduced words of $w$}, denoted $\cR(w)$, is the
undirected graph whose vertex set consists of reduced words for $w$
and edges between reduced words that can be obtained from one another
using a commutation or braid. Letting $\sim_w$ denote commutation
equivalence, one can construct the quotient graph ${\cR(w)/\!\sim}_w$
with vertices given by the equivalence classes of reduced words.  The
\defn{directed braid graph on commutation classes of reduced words of
$w$}, denoted $\cG(w)$, is obtained from $\cR(w)/\!\sim_w$ by
directing the edge between two commutation classes if suitable
representatives from the classes differ by a braid in the direction
$i(i+1)i\rightarrow (i+1)i(i+1)$. See \cref{fig:reducedwordgraph} for
examples.

For $w\in \affineS_n$, the \defn{inversion set} of $w$ is
\begin{equation}\label{eq:inversion.set.def}
\Inv(w)=\left\{(x,y)\in [n] \times \mathbb{Z} : x<y \text{ and } w^{-1}(x)> w^{-1}(y) \right\}.
\end{equation}
Each affine permutation is uniquely
determined by its inversion set, and the length of $w$ can be found using $\ell(w)=|\Inv(w)|$.  The \defn{weak (Bruhat) order} is the partial order on $\affineS_n$ defined by $v\leq w$ whenever $\Inv(v)\subseteq
\Inv(w)$.  The covering relations in the weak order are all of the
form $w \lessdot ws_i$ if $\ell(w) = \ell(ws_i) - 1$ for some $i \in
[n]$, so one could also say the weak order is defined by
\defn{single step inclusion} on inversion sets.  Furthermore, the set of
maximal chains of the interval $[\mathrm{id}, w]$ is in bijection with the set of reduced words for $w$.

Reflection orders for Coxeter groups
were originally studied by Dyer in his Ph.D. thesis \cite{Dyer.1993} and now have been
applied extensively in Coxeter groups and
Kazhdan-Lusztig polynomials. We give a definition in the setting of $\affineS_n$.

\begin{definition}\label{def:reflection order}
Let $w \in \affineS_{n}$.  For each reduced word $\mathbf{i}=i_1
\cdots i_\ell \in \cR(w)$, associate affine permutations
$w^{(j)}=s_{i_{1}}s_{i_{2}}\cdots s_{i_{j-1}}$ for all $1\leq j \leq
\ell$, and define $\rho_{j}=(w^{(j)}(i_{j}), w^{(j)}(i_{j}+1))$ for all
$1\leq j \leq \ell$. The \defn{reflection order} associated to
$\mathbf{i}$ is the total order on $\Inv (w)$ given by
$\rho(\mathbf{i})=(\rho_{1},\dotsc , \rho_{\ell})$.  The set of
\defn{reflection orders} of $w$ is the set of all such
total orders on $\Inv(w)$ associated to reduced words of $w$.
\end{definition}

\begin{example}
If $w=(1,7,2,0) \in \affineS_{4}$, then the reflection order for the
reduced word $0121032$ is
$((0, 1), (0, 2), (0, 3), (2, 3), (1, 3), (0, 7), (2, 7))$
\end{example}

\section{Higher Bruhat orders for permutations in $\S_n$}\label{sec:non-longest}

In this section, we construct the higher Bruhat orders for arbitrary $w\in \S_n$. When $w$ is the longest element in $\S_n$, these reduce to definitions for the classical higher Bruhat orders in \cite{manin_schechtman_1989} and \cite{ziegler_1993}. Throughout this section, assume $1\leq k\leq n$ and $w \in \S_{n}$.

Let $\binom{[n]}{k}$
denote the collection of ordered subsets $[x_1,..., x_k]$ of $k$
distinct elements in $[n]$ listed in increasing order. For each $X = [x_1, \ldots, x_k] \in \binom{[n]}{k}$ and $i\in [k]$, define $X_i=[x_1,\ldots,x_{i-1},x_{i+1},\ldots,x_k]$, and define the \defn{packet} of $X$ to be $P(X) = \{X_1, X_2, \ldots, X_k\}$. This packet has a natural  \defn{lexicographic} (lex) order $(X_k, X_{k-1}, \ldots, X_1)$ and \defn{antilexicographic} (antilex) order $(X_1, X_2, \ldots, X_k)$. A \defn{prefix} of $P(X)$  is a subset of the form $\{X_k, X_{k-1}, \ldots, X_i \}$, and a \defn{suffix} of $P(X)$ is a subset of the form $\{X_i, X_{i-1}, \ldots, X_1 \}$.

In this notation, $[i,j] \in \binom{[n]}{2}$ is an inversion for $w$
whenever $w^{-1}(i)>w^{-1}(j)$.  More generally, we consider the size
$k$-subsets appearing in reverse order for $w$.  

\begin{definition}
For $X=[x_1,\ldots,x_k] \in
\binom{[n]}{k}$, we say $X$ is a ${k}$-\defn{inversion} for
$w$ provided $w^{-1}(x_1)>\cdots >w^{-1}(x_k)$.  The
\defn{$k$-inversion set} of $w$ is defined to 
be
\begin{equation}
{\Inv_k(w)}=\left\{[x_1,\ldots,x_k] \in \binom{[n]}{k}:
w^{-1}(x_1) > \cdots > w^{-1}(x_k) \right\}.
\end{equation}
\end{definition}

 Note that for $k>n$, we have that $\Inv_{k}(w)=\emptyset$, and for $k=0$,  we have that $\Inv_{0}(w)=\{\emptyset \}$.  For $k=1$, every $\{i\}:=[i] \in \binom{[n]
}{1}$ is a 1-inversion of $w \in \S_{n}$. For $k\geq 2$ and $X\in {[n]\choose k+1} \setminus \Inv_{k+1}(w)$ the
intersection $P(X)\cap \Inv_k(w)$ will not be the entire
packet. However, $P(X)\cap \Inv_k(w)$ can be characterized using the following four cases.

\begin{lemma}\label{lem:intersection_size}
     For any $X \in \binom{[n]}{k+1}$, the intersection $P(X)\cap
\Inv_k(w)$ is one of the
following: $\emptyset$, 
$\{X_i\}$ for some $i\in [k+1]$, $\{X_i,X_{i+1}\}$ for some
$i\in [k]$, or the entire packet $P(X)$. 
\end{lemma}

\cref{lem:intersection_size} motivates another definition. For $X=[x_1,\ldots, x_k] \in \binom{[n]}{k}$, we say $X$ is a \defn{$k$-quasi-inversion} for $w$ if all but one of the pairs in $\{[x_i, x_j] : 1 \le i < j \le k\}$ are 2-inversions of $w$. When this occurs, $P(X) \cap \Inv_{k-1}(w)=\{X_i,X_{i+1}\}$ for some $i\in [k-1]$.
\begin{example}\label{ex:k.inv}
    For $w=(6,4,5,2,3,1)\in \S_6$, we see $[1,2,4,6]$ is a $4$-inversion, $[1,2,3,4]$ is a $4$-quasi-inversion, and $[2,3,4,5]$ is neither an inversion nor a quasi-inversion. 
\end{example}

A $2$-quasi-inversion is sometimes referred to as a co-inversion of $w$. Furthermore, permutations in the interval $[\mathrm{id},w]$ in the weak order can be identified with the linear extensions of a poset on $[n]$ with
relations $i<j$ if $[i,j]$ is a quasi-inversion for $w$.  For larger $k$, we will
construct a similar partial order based on quasi-inversions.

\begin{definition}\label{def:permanent.poset.}
The \defn{permanent poset} ${\P_w(n,k)}$ is the poset on $\Inv_{k}(w)$ with order
relation $\leq_{\P_w(n,k)}$ given by the transitive closure of the \defn{quasi-inversion
relations}: if $X \in \binom{[n]}{k+1}$ is a quasi-inversion for $w$
such that $P(X)\cap \Inv_{k}(w) = \{X_i,X_{i+1}\}$, then $X_{i}<_{\P_w(n,k)} X_{i+1}$ whenever $(k-i)$ is odd, and
$X_{i+1}<_{\P_w(n,k)} X_{i}$ whenever $(k-i)$ is even.
\end{definition}

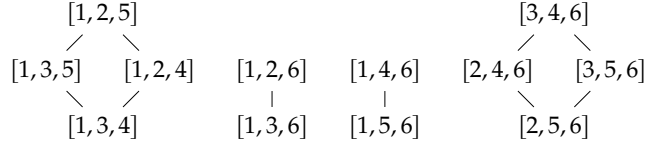
\begin{figure}[t]
\centering 
\scalebox{0.75}{
    \begin{tikzpicture}
        \node (134) at (0,0) {$[1,3,4]$};
        \node (135) at (-1,1) {$[1,3,5]$};
        \node (124) at (1,1) {$[1,2,4]$};
        \node (125) at (0,2) {$[1,2,5]$};
        \node (136) at (3,0) {$[1,3,6]$};
        \node (126) at (3,1) {$[1,2,6]$};
        \node (156) at (5,0) {$[1,5,6]$};
        \node (146) at (5,1) {$[1,4,6]$};
        \node (256) at (8,0) {$[2,5,6]$};
        \node (246) at (7,1) {$[2,4,6]$};
        \node (356) at (9,1) {$[3,5,6]$};
        \node (346) at (8,2) {$[3,4,6]$};
        \path (134) edge node {} (135);
        \path (134) edge node {} (124);
        \path (135) edge node {} (125);
        \path (124) edge node {} (125);
        \path (136) edge node {} (126);
        \path (156) edge node {} (146);
        \path (256) edge node {} (246);
        \path (246) edge node {} (346);
        \path (256) edge node {} (356);
        \path (356) edge node {} (346);
    \end{tikzpicture}}
    \caption{The Hasse diagram of $\P_w(6,3)$ for $w=(6,4,5,2,3,1)\in \S_6$.}
    \label{fig:permanent example}
\end{figure}

To give some further intuition on the relations, consider the reflection orders for $w$ in \cref{def:reflection order}. One can show that for any quasi-inversion $[x,y,z]$, if $P([x,y,z])\cap \Inv_2(w)=\{[x,y],[x,z]\}$ (resp. $\{[x,z],[y,z]\}$), then $[x,y]$ (resp. $[y,z]$) must appear before $[x,z]$. These aligns with the relations in $\P_w(n,2)$. See \cref{fig:permanent example} for an example of $\P_w(n,3)$.

Reflexivity and transitivity of $\leq_{\P_w(n,k)}$ follow by construction. To establish antisymmetry and show that $\P_w(n,k)$ actually is a poset, we construct functions that are monotone on $\P_w(n,k)$. These functions also allow us to describe the restriction of $\P_w(n,k)$ to $P(X)$ for any $X\in \Inv_{k+1}(w)$.

\begin{lemma}\label{lem:poset.packet.antichain}
The permanent poset $\P_{w}(n,k)$ is a poset.  Furthermore, if $X\in
\Inv_{k+1}(w)$, then the elements in $P(X)$ form an antichain in
$\P_w(n,k)$.
\end{lemma}

We now consider admissible orders. The primary modification from the
classical admissible orders is a separate condition for
quasi-inversions, which is encoded in $\P_w(n,k)$. 

\begin{definition}\label{def:admissible.BCL}
A total
order $\rho$ of $\Inv_k(w)$ is a ${k}$-\defn{admissible
order} for $w$ if $\rho$ is a linear extension of $\P_w(n,k)$, and for every $X\in \Inv_{k+1}(w)$, the restriction $\rho|_{P(X)}$ is the lex or antilex order on $P(X)$.  Let ${\A_w(n,k)}$ be the \defn{set of $k$-admissible orders}
    of $\Inv_k(w)$.  
    Each admissible order has a \defn{reversal set} defined by
    \begin{equation}
        {\Rev_{n,k,w}}(\rho)=   \{X \in \Inv_{k+1}(w):  \rho|_{P(X)}=(X_1,X_2,\ldots,X_{k+1})\}.
    \end{equation}

\end{definition}

One can show that the admissible orders in $\A_{w}(n,1)$ are
exactly the permutations in the interval $[\mathrm{id},w]$ in the weak
order on $\S_{n}$ if one equates each singleton set $\{i \}$ with the
value $i$.  Additionally, the admissible orders in $\A_{w}(n,2)$ are the
reflection orders for $w$. Examples of admissible orders for $k=3$ are shown in \cref{table:admisible}.

\begin{table}[t]
\centering 
\scriptsize 
    \begin{tabular}{|c|c|}  \hline 
    Admissible order & Reversal Set \\ \hline 
        $\rho^1 = (134, 124, 135, 136, 125, 126, 156, 256, 146, 246, 356, 346)$ & $\emptyset$\\ \hline 
        $\rho^2 = (134, 256, 135, 136, 156, 356, 124, 126, 146, 246, 346, 125)$ & $\{1256\}$\\ \hline 
        $\rho^3 = (134, 256, 135, 136, 156, 356, 246, 146, 126, 124, 346, 125)$ & $\{1246, 1256\}$\\\hline 
        $\rho^4 = (134, 256, 356, 156, 136, 135, 124, 126, 146, 246, 346, 125)$ & $\{1256, 1356\}$\\ \hline 
        $\rho^5 = (134, 256, 356, 156, 136, 135, 246, 146, 126, 124, 346, 125)$ & $\{1246, 1256, 1356\}$\\\hline 
        $\rho^6 = (256, 246, 356, 156, 346, 146, 136, 134, 126, 124, 135, 125)$ & $\Inv_4(w)$ \\\hline 
    \end{tabular}
    \caption{Six $3$-admissible orders for the permutation $w=(6,4,5,2,3,1)\in \S_6$ from \cref{fig:permanent example}.  Square brackets and internal commas have been suppressed for clarity.}
\label{table:admisible}
\end{table}

We next define two operations on admissible orders. Elements $X,Y\in \Inv_k(w)$ \defn{commute with respect to}
${w}$ if they are incomparable in $\P_{w}(n,k)$ and do not belong to a common packet.  Two total orders of $\Inv_{k}(w)$ are \defn{commutation equivalent} if they can be obtained from one another using commutations on adjacent elements that commute with respect to $w$. We will extend $\sim_{w}$ from reduced words to admissible orders to denote this equivalence relation, and the commutation class of $\rho\in \A_w(n,k)$
will be denoted ${[\rho]}$. 

For any $\rho\in \A_w(n,k)$ and $X\in \Inv_{k+1}(w)$, the packet
$P(X)$ is \defn{flippable} in $\rho$ if $P(X)$ forms a saturated chain
in the total order $\rho$.  If $\rho|_{P(X)}$ is  the lex (resp. antilex) order on $P(X)$, then a
\defn{lex-to-antilex packet flip} (resp. \defn{antilex-to-lex packet
flip}) at $P(X)$ is the operation that reverses the order of $P(X)$ in
$\rho$.  Similarly, the packet $P(X)$ is \defn{flippable} for the
equivalence class $[\rho]\in \A_w(n,k)/\!\sim_{w}$ if it is flippable for some
representative in $[\rho]$.

\begin{example}\label{ex:flips}
    In the admissible order $\rho^2$ from \cref{table:admisible},
each of $P(1246)$ and $P(1356)$ form a saturated chain in
lex order, so each packet is lex-to-antilex flippable. Flipping these
respectively results in $\rho^3$ and $\rho^4$. Additionally, $\rho^2$
is commutation equivalent to
\begin{equation}
    \sigma = (134, 124, 135, 136, 256,156, 126, 125, 146, 246, 356, 346).
\end{equation}
Observe that $P(1256)$ is antilex-to-lex flippable for $\sigma$,
and flipping this results in $\rho^1$.  
\end{example}

\begin{lemma}\label{lem:reversal_sets}
Let $\rho\in \A_w(n,k)$.
\begin{enumerate}[label=(\alph*)]
    \item The reversal set $\Rev_{n,k,w}(\rho) \subseteq
\Inv_{k+1}(w)$ is a (lower) order ideal of $\P_w(n,k+1)$, and its
intersection with packet $P(X)$ for $X\in \binom{[n]}{k+2}$ is a prefix or suffix of $P(X)$.
    \item If $\sigma \sim_{w} \rho$, then $\sigma \in
\A_w(n,k)$ and $\Rev_{n,k,w}(\rho) = \Rev_{n,k,w}(\sigma)$.
    \item If $\sigma$ is obtained from $\rho$ by a packet flip at $P(X)$, then we have that $\sigma\in \A_w(n,k)$ and $\Rev_{n,k,w}(\sigma)\bigtriangleup \Rev_{n,k,w}(\rho)=\{X\}$, where $\bigtriangleup$ denotes symmetric difference. 
\end{enumerate}
\end{lemma}

\cref{lem:reversal_sets} shows that
$\A_w(n,k)$ is closed under commutations and
packet flips. With this in mind, we now define the higher Bruhat orders for any $w\in \S_n$.

\begin{definition}\label{def:higher.B.Sn}
The
${k^{th}}$\defn{ higher Bruhat order for $w$,} denoted 
${\B_w(n,k)}$,  is the partial order on $\A_w(n,k)/\!\!\sim_{w}$ where $[\rho]\leq [\sigma]$ if $\sigma$ can be obtained from
$\rho$ by some sequence of commutations and lex-to-antilex packet
flips.
\end{definition}

\begin{figure}[t]
    \centering 
    \scalebox{0.75}{
    \begin{tikzpicture}
        \node (rho1) at (0,0) {$[\rho^1]$};
        \node (rho2) at (0,1) {$[\rho^2]$};
        \node (rho4) at (1,2) {$[\rho^4]$};
        \node (rho3) at (-1,2) {$[\rho^3]$};
        \node (rho5) at (0,3) {$[\rho^5]$};
        \node (rho6) at (0,4) {$[\rho^6]$};
        \path (rho1) edge node {} (rho2);
        \path (rho2) edge node {} (rho3);
        \path (rho2) edge node {} (rho4);
        \path (rho3) edge node {} (rho5);
        \path (rho4) edge node {} (rho5);
        \path (rho5) edge node {} (rho6);
        \node (1256) at (5,1) {$[1,2,5,6]$};
        \node (1246) at (4,2) {$[1,2,4,6]$};
        \node (1356) at (6,2) {$[1,3,5,6]$};
        \node (1346) at (5,3) {$[1,3,4,6]$};
        \path (1256) edge node {} (1246);
        \path (1256) edge node {} (1356);
        \path (1246) edge node {} (1346);
        \path (1356) edge node {} (1346);
        \node (r1) at (12.5,0) {$\emptyset$};
        \node (r2) at (12.5,1) {$\{[1,2,5,6]\}$};
        \node (r4) at (10.5,2) {$\{[1,2,4,6],[1,2,5,6]\}$};
        \node (r3) at (14.5,2) {$\{[1,2,5,6],[1,3,5,6]\}$};
        \node (r5) at (12.5,3) {$\{[1,2,4,6],[1,2,5,6],[1,3,5,6]\}$};
        \node (r6) at (12.5,4) {$\Inv_4(w)$};
        \path (r1) edge node {} (r2);
        \path (r2) edge node {} (r3);
        \path (r2) edge node {} (r4);
        \path (r3) edge node {} (r5);
        \path (r4) edge node {} (r5);
        \path (r5) edge node {} (r6);
    \end{tikzpicture}}
    \caption{The Hasse diagrams for $\B_w(6,3)$ (left), $\P_w(6,4)$ (middle), and $\C_w(6,4)$ (right) for the permutation $w=(6,4,5,2,3,1)$ from \cref{table:admisible}. }
    \label{fig:BnkCnk}
\end{figure}
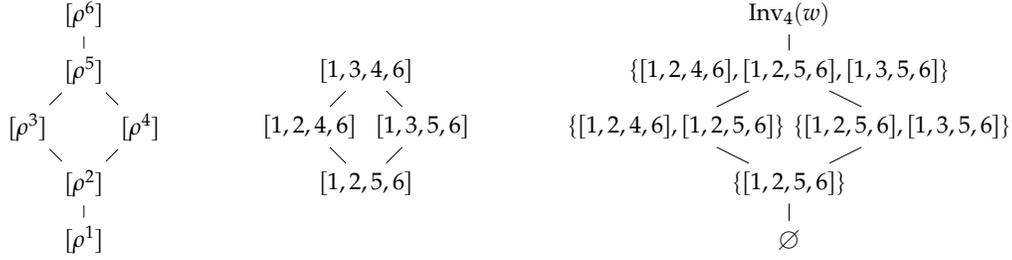

From the definition, it is straightforward to show $B_{w}(n,1)$ is isomorphic to the weak order
interval $[\mathrm{id},w]$. We will see later from \cref{thm:nonlongest} that the Hasse diagram of $\B_{w}(n,2)$ can be viewed as the graph $\cG(w)$. For an example with $k=3$, one can use a computer to verify that there are 1228 $3$-admissible orders for $w=(6,4,5,2,3,1)$, and the six admissible orders in \cref{table:admisible} give a
representative from each of the six commutation classes in $\A_w(6,3)/\!\sim_w$. The poset $\B_w(6,3)$ is shown on the left in \cref{fig:BnkCnk}.

\cref{lem:reversal_sets}(a) also describes necessary conditions for reversal sets. Based on this, we extend the notion of consistent sets from \cite{ziegler_1993}.  See \cref{fig:BnkCnk} for examples.

\begin{definition}\label{def:consistent.Sn}
    A subset $R\subseteq \Inv_{k}(w)$ is \defn{consistent with respect to ${w}$} if $R$ is an order ideal of $\P_w(n,k)$ that satisfies the \defn{Manin-Schechtman-Ziegler (MSZ) Condition}: \[\text{for any $X\in \Inv_{k+1}(w)$, the intersection $P(X)\cap R$ is a prefix or suffix of $P(X)$}.\]
Define ${\C_w(n,k)}$ to be the poset on consistent subsets
of $\Inv_k(w)$ with partial order generated by single step
inclusion. 
\end{definition}

We will conclude this section by outlining our approach for proving \cref{thm:nonlongest}. By \cref{lem:reversal_sets}, 
$\Rev_{n,k,w}:\A_w(n,k)\to \C_w(n,k+1)$ is well-defined and descends to the quotient $\B_w(n,k)$. To prove this is a bijection, we must show each 
$R\in \C_w(n,k+1)$ is the reversal set for some admissible order in
$\A_w(n,k)$, and each $[\rho]\in \B_w(n,k)$ is uniquely determined by its reversal set. We do this by considering a directed
graph $G_R$ on $\Inv_k(w)$ constructed from 
$R\in \C_w(n,k+1)$. This graph has directed edges
\begin{itemize}
\item $X \to Y$ if $X<_{{\P_{w}(n,k)}}Y$ is a quasi-inversion
relation,
\item $X_i\color{red}\to\color{black} X_{i+1}$ for all $1\leq i\leq  k$
if $X\in R$, and
\item $X_{i+1}\color{blue}\to\color{black} X_i$ for all $1\leq i\leq  k$ if $X\in\Inv_{k+1}(w)\setminus R$. \end{itemize}

    \begin{figure}[b]
\centering \scalebox{0.75}{
    \begin{tikzpicture}
        \node (134) at (0,0) {$[1,3,4]$};
        \node (135) at (-1,1) {$[1,3,5]$};
        \node (124) at (1,1) {$[1,2,4]$};
        \node (125) at (0,2) {$[1,2,5]$};
        \node (136) at (3,0) {$[1,3,6]$};
        \node (126) at (3,1) {$[1,2,6]$};
        \node (156) at (5,0) {$[1,5,6]$};
        \node (146) at (5,1) {$[1,4,6]$};
        \node (256) at (8,0) {$[2,5,6]$};
        \node (246) at (7,1) {$[2,4,6]$};
        \node (356) at (9,1) {$[3,5,6]$};
        \node (346) at (8,2) {$[3,4,6]$};
        \draw[->] (134) edge node {} (135);
        \draw[->] (134) edge node {} (124);
        \draw[->] (135) edge node {} (125);
        \draw[->] (124) edge node {} (125);
        \draw[->] (136) edge node {} (126);
        \draw[->] (156) edge node {} (146);
        \draw[->] (256) edge node {} (246);
        \draw[->] (246) edge node {} (346);
        \draw[->] (256) edge node {} (356);
        \draw[->] (356) edge node {} (346);
        \draw[<-,color=red] (125) edge node {} (126);
        \draw[<-,color=red] (126) edge node {} (156);
        \draw[<-,color=red] (156) edge node {} (256);
        \draw[<-,color=red] (135) edge node {} (136);
        \draw[<-,color=red] (136) edge node {} (156);
        \draw[<-,color=red] (156) edge node {} (356);
        \draw[->,color=blue] (124) edge node {} (126);
        \draw[->,color=blue] (126) edge node {} (146);
        \draw[->,color=blue] (146) edge node {} (246);
        \draw[->,color=blue] (134) edge node {} (136);
        \draw[->,color=blue] (136) edge node {} (146);
        \draw[->,color=blue] (146) edge node {} (346);
    \end{tikzpicture}}
    \caption{The graph $G_R$ for $w=(6,4,5,2,3,1)$ and $R=\{[1,2,5,6],[1,3,5,6]\}\in\C_w(6,4)$.}
    \label{fig:Gr example}
\end{figure}
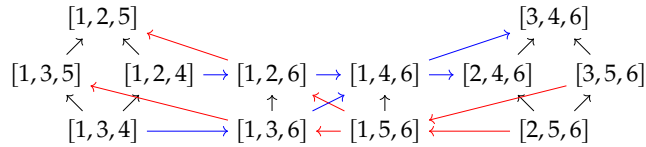

An example is shown in \cref{fig:Gr example}. Observe that the structure of $G_R$ can be complex.  Through a technical induction argument on $n$ and $k$ that decomposes $G_R$ into well-behaved subgraphs, we establish the following structural result.

\begin{lemma}\label{lem:GRacyclic}
For any $R\in \C_w(n,k+1)$, the directed graph
$G_R$ is acyclic.
\end{lemma}

As $G_R$ is acyclic, it induces a partial order $\leq_R$ on $\Inv_k(w)$. Letting $\mathscr{L}(\Inv_k(w),\leq_R)$ denote the set of linear extensions of $(\Inv_k(w),\leq_R)$, we establish the following result.

\begin{theorem}\label{thm:Bnk_iso_Cnk}
    The map $\Rev_{n,k,w}: \B_w(n,k) \to
\C_w(n,k+1)$ is a poset isomorphism with inverse given
by
$\Rev_{n,k,w}^{-1}(R)=\mathscr{L}(\Inv_k(w),\leq_R).$
\end{theorem}

To show that $\C_w(n,k)$ is ranked with unique minimal element $\emptyset$ and maximal element $\Inv_{k}(w)$, it suffices to show that a $k$-inversion can be removed or added whenever $R\in \C_w(n,k)$ is not $\emptyset$ or $\Inv_{k}(w)$, respectively. We do this by defining the \defn{suffix set} for $R$ as $S(R) = \{Y \in \Inv_{k+1}(w): Y_1\in R\}$. After showing $S(R)\in \C_w(n,k+1)$, we use $G_{S(R)}$ to establish the following result.

\begin{corollary}\label{thm:Cnk_min_max}
The poset $\C_w(n,k)$
is a ranked poset with unique minimal element $\emptyset$ and a unique
maximal element $\Inv_{k}(w)$.  The rank of $R \in \C_w(n,k)$ is
$|R|$.
\end{corollary}

For \cref{thm:nonlongest}, it remains to show the natural bijection
between maximal chains of $\C_w(n,k+1)$ and $\A_w(n,k+1)$ is given by
using the single step inclusions 
\begin{equation} 
    (\emptyset\lessdot R_1 \lessdot R_2 \lessdot \dots \lessdot R_{|\Inv_{k}(w)|})\mapsto (R_1,R_2\setminus R_1,\ldots,R_{|\Inv_{k}(w)|}\setminus R_{|\Inv_{k}(w)|-1}).
\end{equation}
 While proving injectivity is straightforward, proving surjectivity requires a technical argument that utilizes consequences of \cref{thm:Bnk_iso_Cnk,thm:Cnk_min_max}.

\section{Higher Bruhat orders for permutations in $\affineS_n$}\label{sec:affine}

In this section, we will generalize the definitions from
\Cref{sec:non-longest} to elements of the affine symmetric
group. Throughout, assume $2\leq n$,\ $1\leq k\leq n$, and $w\in
\affineS_n$.

In order to generalize the notion of $k$-inversions to affine
permutations accounting for $n$-periodicity, let $\binom{\Z }{k}$ denote
the $k$-subsets of $\Z$.  Define an equivalence relation $\sim_n$ on
the elements of $\binom{\Z }{k}$ where $\{x_1, x_2, \ldots, x_k \}
\sim_n \{y_1, y_2, \ldots, y_k \}$ if $k>1$ and there exists an
integer $m$ such that $\{x_1, x_2, \ldots, x_k \} = \{y_1+mn, y_2+mn,
\ldots, y_k+mn\}$. If $k=1$, then $\{x \} \sim_n \{y \}$ if and only
if $x = y$.  Denote the affine equivalence classes of size $k$ subsets of
$\Z$ with distinct elements modulo $n$ listed in increasing order as
\begin{equation}
    \binom{\mathbb{Z}}{k}_n = \bigl\{ [x_1, \ldots, x_k] : x_i
\not\equiv x_j \pmod{n}\text{ for $1 \le i < j \le k$} \bigr\}.
\end{equation}
For each $X\in  \binom{\mathbb{Z}}{k}_n$, the packet of $X$ is defined similarly as in \cref{sec:non-longest}.

Adapting the definition from \cref{sec:non-longest}, $X=[x_1,\ldots,x_k] \in
\binom{\Z}{k}_{n}$ is a ${k}$-\defn{inversion} for $w$
provided $w^{-1}(x_1)>\cdots >w^{-1}(x_k)$.  We then define the \defn{$k$-inversion
set} of $w$ to be
\begin{equation}
    \Inv_k(w)=\left\{[x_1,\ldots,x_k]\in \binom{\mathbb{Z}}{k}_n:
w^{-1}(x_1)>\ldots >w^{-1}(x_k) \right\}.
\end{equation}

Using the notation of affine equivalences classes in $\binom{\mathbb{Z}}{k}_n$, the following statement is a characterization
of $2$-inversion sets for affine permutations. The proof follows from \cite[Prop
2.1]{Barkley.Speyer.2024} and \cite[Lemma 4.1(d)]{Dyer.2019}.

\begin{theorem}\label{lem:inversion.sets} Let $R\subseteq \binom{\mathbb{Z}}{2}_n$.
Then, $R$ is the inversion set for some affine permutation in
$\affineS_{n}$ if and only if for all $[x,y,z]\in \binom{\mathbb{Z}}{3}_n$, we have
\begin{itemize}
\item $[x,z] \in R$ implies  $[x,y] \in R$ or $[y,z] \in R$, 
\item $[x,y] \in R$ and $[y,z] \in R$ implies  $[x,z] \in R$, and
\item $[x,y] \in R$ implies $[x,y-en]
\in R$ for all $e \in \Z_{\geq 0}$ such that $ x\leq y-en\leq y$.
\end{itemize}
\end{theorem}

The reflection orders for $w\in \affineS_n$ have the property that each ordered prefix must also be
a reflection order since each ordered prefix of a reduced word is still a
reduced word.  Therefore, one can use \Cref{lem:inversion.sets} to characterize reflection orders.

\begin{corollary}\label{cor:characterize.ref.order} Let $w \in
\affineS_{n}$.  A total order $\rho = (\rho_{1},\rho_{2},\dotsc ,
\rho_{\ell})$ of $\Inv_2(w)$ is a reflection order if and only if
\begin{itemize}
\item for each  $[x,y,z] \in \binom{\mathbb{Z}}{3}_n$ the total order
$\rho$ restricted to $P([x,y,z])\cap \Inv_2(w)$ is a prefix of
$P([x,y,z])$ ordered in lex order or a suffix of $P([x,y,z])$ ordered in antilex order, and
\item for each pair $[x,y],[x,y+n] \in \Inv_2(w)$, the pair $[x,y]$ appears before
$[x,y+n]$ in $\rho$.
\end{itemize}
\end{corollary}

The definition of quasi-inversion in \cref{sec:non-longest} carries over similarly from $\S_n$ to $\affineS_n$. As in $\S_n$, we define a permanent poset for $w \in
\affineS_{n}$. This is motivated by the weak order on
$\affineS_{n}$, computer experimentation,  \Cref{lem:inversion.sets}, and \Cref{cor:characterize.ref.order}.

\begin{definition}\label{def:congruence.poset}
    Let $v_i^{(n, k)}$ denote the vector $(0, \ldots, 0, n, \ldots, n)$
consisting of $i$ copies of 0 followed by $(k-i)$ copies of $n$. The \defn{congruence poset} $\CR(n,k)$ is the poset on $\binom{\Z}{k}_n$ generated by the
\defn{congruence relations} for all $X \in \binom{\Z}{k}_n$ and $0 \leq 
i < k$: $X \leq_{\CR(n,k)} X + v_i^{(n,k)}$ if $k-i$ is odd, and $X + v_i^{(n,k)} \leq_{\CR(n,k)} X$ if $k-i$ is even.
The \defn{permanent poset}
${\P_w(n,k)}$ is the set $\Inv_{k}(w)$ with order
relations given by the transitive closure of the quasi-inversion
relations of \cref{def:permanent.poset.} and the congruence relations
of $\CR(n,k)$ restricted to $\Inv_{k}(w)$.  
\end{definition}

When $k = 1$, we define $\A_w(n,1)$ in such a way that $1$-admissible
orders for $w$ are in bijection with the interval $[\mathrm{id},w]$ in
weak order.  A total order $\rho$ on $\Inv_1(w)\cong \Z$ is a 1-\defn{admissible
order} if $\rho$ is a linear extension of $\P_w(n,1)$ with a
well-defined finite reversal set in $\binom{\Z }{2}_{n}$. So
$[z_{1},z_{2}] \in \Rev(\rho)$ implies every pair $\{z_1+mn,z_2+mn\}$
for all $m\in \mathbb{Z}$ appears in antilex order in $\rho$.  With
these definitions in mind, one can also show that
$\B_w(n,1)\cong\C_w(n,2)$ is isomorphic to $[\mathrm{id},w]$ in the
weak order on $\affineS_n$.

The ${k}$-{admissible orders} $\A_w(n,k)$ for $k\geq 2$ are defined as
linear extensions of $\P_w(n,k)$ that satisfy the MSZ condition as in
\cref{def:admissible.BCL}, but using the affine notion of
$\P_w(n,k)$. Reversal sets, commutations, and packet flips for
admissible orders are all defined as in \cref{sec:non-longest}, using
$\binom{\mathbb{Z}}{k}_n$ in place of $\binom{[n]}{k}$. Definitions
for $\B_w(n,k)$ as a poset on $\A_w(n,k)/\!\sim_w$ and $\C_w(n,k)$ as
a poset on consistent subsets of $\Inv_k(w)$ carry over mutatis
mutandis. The following result establishes the connections between
$\A_w(n,2)$, reduced words, and reflection orders.

\begin{lemma}\label{lem:Bwn2_and_graph}
    For any $w\in \affineS_n$, the admissible orders in $\A_w(n,2)$ are in natural bijection with reflection orders and reduced words of $w$. Under this bijection, two reduced words
     \begin{enumerate}[label=(\alph*)]
            \item are commutation equivalent if and only if their corresponding two admissible orders are commutation equivalent, and
            \item differ by a braid $s_is_{i+1}s_i\to s_{i+1}s_is_{i+1}$ if and only if their corresponding admissible orders differ by a lex-to-antilex  packet flip. 
    \end{enumerate}
\end{lemma}

\begin{example}
Let $w=(-3, -2, 8, 7)\in \affineS_4$.  
An
admissible order in $\A_w(4,2)$ corresponding to the reduced word $232124134$ 
is 
\begin{equation}
    \rho = ([2,3],[2,4],[3,4],[1,4],[1,3],[2,8],[2,7],[1,8],[1,7])
\end{equation}
 with reversal set
$\Rev_{n,k,w}(\rho) = \{[1,3,4], [2,7,8], [1,7,8]\}$. Observe that
$\Rev_{n,k,w}(\rho)$ is an order ideal of $\P_w(4,3)$ and is a consistent
subset of $\Inv_3(w)$. Applying the braid move $232124134\to 323124134$ corresponds to a lex-to-antilex packet flip at $P([2,3,4])$.
\end{example}

One can show generalizations of \cref{lem:poset.packet.antichain,lem:reversal_sets} hold for $w\in \affineS_n$. However, for $R\in \C_w(n,k)$, we must extend the definition of the graph $G_R$ by including edges $X\to Y$ when $X<_{\P_w(n,k)} Y$ is a congruence relation. The key obstruction to proving a complete analog of \cref{thm:nonlongest} is our proof that $G_R$ is acyclic, which cannot be directly adapted to affine permutations due to the congruence relations. However, we conjecture that an analog of \cref{lem:GRacyclic} does hold. As noted in the Introduction, we have computationally verified this conjecture for  $1\leq k\leq 8$ and affine permutations $w\in \affineS_n$ up to certain lengths.

\begin{conjecture}\label{conj:affineG_R}
For any $R\in \C_w(n,k)$, the graph $G_R$ is acyclic. 
\end{conjecture}

Through a technical argument, we are able to prove the special case of $k=3$. Combining this with \cref{lem:Bwn2_and_graph} and generalizations of results in \cref{sec:non-longest}, we establish \cref{thm:affine}. We note that if \cref{conj:affineG_R} can be resolved for general $k$, then this may lead to a full generalization of \cref{thm:nonlongest} to $w\in \affineS_n$.

\acknowledgements{
We would like to thank Ben Elias for suggesting the problems in the
introduction to us and many helpful conversations along the way.   
}

\printbibliography

\end{document}